\documentclass{amsart}

\usepackage{amssymb,amsmath,amsthm,color}

\usepackage[bookmarksopen=true,bookmarksnumbered=true, bookmarksopenlevel=2, colorlinks=true,linkcolor=mblue,
citecolor=mblue,urlcolor=mblue, pdfstartview=FitH]{hyperref}

\definecolor{mblue}{rgb}{0,0,.8}

\newtheorem{theorem}{Theorem}
\newtheorem{lemma}{Lemma}
\newtheorem{prop}{Proposition}

\newtheorem{conjecture}{Conjecture}

\newtheorem{rmk}{Remark}

\newcommand{\Z}{\mathbb Z}
\newcommand{\Q}{\mathbb Q}
\newcommand{\F}{\mathbb F}

\newcommand{\C}{\mathbb C}
\newcommand{\chibar}{\overline{\chi}}
\newcommand{\sigmabar}{\overline{\sigma}}
\newcommand{\ebar}{\overline{\varepsilon}}
\newcommand{\p}{{\mathfrak p}}

\def\GL{\mathop{\mathrm{GL}}\nolimits}

\title{Mod $pq$ Galois representations and Serre's conjecture}

\author{Chandrashekhar Khare \& Ian Kiming}

\address{C.K.: Department of Mathematics, University of Utah, 155 S 1400 E, Salt Lake City,
UT 84112-0090, USA.}
\email{\href{mailto:shekhar@math.utah.edu}{shekhar@math.utah.edu}}
\address{I.K.: Department of Mathematics, University of Copenhagen, Universitetsparken 5, DK-2100 Copenhagen \O ,
Denmark.}
\email{\href{mailto:kiming@math.ku.dk}{kiming@math.ku.dk}}

\begin{document}

\begin{abstract} We consider linear representations of the
Galois groups of number fields in 2 different characteristics
and examine conditions under which they arise simultaneously
from a motive.
\end{abstract}

\maketitle

\section{Introduction.}\label{intro}
Motives and automorphic forms of arithmetic type give rise to Galois representations that occur in {\it compatible families}.
These compatible families are of $p$-adic representations (see \cite{S}) with $p$ varying. By reducing such a family mod $p$ one
obtains compatible families of mod $p$ representations (see \cite{K-Pre}). While the representations that occur in such a
$p$-adic or mod $p$ family are strongly correlated, in a sense each member of the family reveals a new face of the motive. In
recent celebrated work of Wiles playing off a pair of Galois representations in different characteristics has been crucial.

In this paper we investigate when a pair of mod $p$ and mod $q$ representations of the absolute Galois group of a
number field $K$ simultaneously arises from an {\it automorphic motive}: we do this in the 1-dimensional (Section 2)
and 2-dimensional (Section 3: this time assuming $K={\mathbb Q}$) cases. In Section 3 we formulate a mod $pq$ version
of Serre's conjecture refining in part the question of Barry Mazur and Ken Ribet in \cite{Stein} and some of the
considerations in \cite{K-MRL}.

\subsection{Notation.}\label{notation}
Throughout the paper, $p$ and $q$ denote 2 distinct, fixed prime
numbers.
\smallskip

We fix the following notation:
$$\begin{array}{lll} K & : & \mbox{a finite extension of } \Q \\
E & : & \mbox{group of global units in } K\\
G_K & : & \mbox{the absolute galois group of } K\\
C_K & : & \mbox{the idele class group of } K\\
\Sigma & : & \mbox{the set of embeddings } K\hookrightarrow
\overline\Q \end{array}$$

By `Hecke character' of $K$ we shall mean a grossencharacter of
type $A_0$ on $C_K$.
\smallskip

We fix embeddings $\iota_{\infty} \colon \overline\Q \hookrightarrow \C$, and for every prime $\ell$ an embedding $\iota_{\ell}
\colon \overline\Q \hookrightarrow \overline\Q_{\ell}$, and denote by $v_{\ell}$ the extension of the $\ell$-adic valuation on
$\Q_{\ell}$ to $\overline\Q_{\ell}$. We also fix an isomorphism of the residue field of $\overline\Q_{\ell}$ with
$\overline\F_{\ell}$.

Thus, an element $\sigma \in \Sigma$ gives rise to an embedding
$K\hookrightarrow \overline\Q_{\ell}$ via composition with
$\iota_{\ell}$. Abusing notation, with denote again this
embedding by $\sigma$.

If $\lambda$ is a finite prime of $K$ of residue characteristic
$\ell$ we denote by:
$$\begin{array}{lll} {\mathfrak p}_{\lambda} & : & \mbox{maximal
ideal of ring of integers in $K_{\lambda}$} \\
k_{\lambda} & : & \mbox{residue field of } \lambda \\
U_{\lambda} & : & \mbox{corresponding unit group}
\end{array}$$ and introduce the
subset $\Sigma(\lambda) \subseteq \Sigma$:
$$\begin{array}{lll} \Sigma(\lambda) & := &
\{ \sigma \in \Sigma \mid ~ v_{\ell} \circ \sigma \cong \lambda \}. \end{array}$$Thus, if $\sigma \in \Sigma(\lambda)$ we have a
canonical extension to an embedding $\sigma_{\lambda} \colon K_{\lambda} \hookrightarrow \overline\Q_{\ell}$ inducing an
embedding:$$\sigmabar_{\lambda} \colon k_{\lambda} \hookrightarrow \overline\F_{\ell} .$$

Notice that the sets $\Sigma(\lambda)$ form a partition of
$\Sigma$ when $\lambda$ runs over the prime divisors of $\ell$.
\smallskip

Since we have fixed a prime of $\overline\Q$ above $\lambda$, we
can speak of the inertia group $I_{\lambda} \le G_K$ at $\lambda$.
We denote by$$\theta_{\lambda} \colon I_{\lambda} \longrightarrow
k_{\lambda}^{\times}$$the character (factoring through tame
inertia) that gives the action of $\mathrm{Gal}(\overline
K_{\lambda}/K_{\lambda})$ on a $(\# k_{\lambda} -1)$'st root of a
uniformizer of $K_{\lambda}$. Thus, if $\# k_{\lambda} = \ell^n$
the fundamental characters of level $n$ over $K_{\lambda}$ are
the $n$ characters:$$\theta_{\lambda}^{\ell^i}, \qquad \mbox{for}
\quad i = 0, \ldots ,n-1.$$

We denote by $\tilde\theta_{\lambda}$ the Teichm\"{u}ller lift of $\theta_{\lambda}$, and --- as we have fixed an embedding
$\overline\Q \hookrightarrow \overline K_{\lambda}$ --- can regard it as taking values in either $\overline\Q^{\times}$, or in
$\overline K_{\lambda}^{\times}$.
\smallskip

In order to reduce the amount of notation, we allow ourselves to
use class field theory implicitly. Thus, if $\chi : G_K
\rightarrow \overline\F_{\ell}^{\times}$ is a homomorphism we may
find occasion to view $\chi$ as a homomorphism on $C_K$.

Generally, if $\chi$ is a Hecke character of $K$ and if $\ell$ is
a prime number, there is attached to the pair $(\chi,\ell)$ an
$\ell$-adic galois representation $\chi_{\ell} : G_K \rightarrow
\overline\Q_{\ell}^{\times}$ satisfying appropriate reciprocity
laws (cf. for instance \cite{S1}, \cite{S}, or the original
construction \cite{weil}). The character $\chi_{\ell}$ lands in
the group of units. Reducing `mod $\ell$' we obtain a
representation $\chibar_{\ell}$:$$\chibar_{\ell} : G_K
\longrightarrow \overline\F_{\ell}^{\times}.$$In particular, we
use these notations $\varepsilon_{\ell}$ and
$\bar\varepsilon_{\ell}$ if $\varepsilon$ is just a finite-order
complex character of $C_K$.
\smallskip

Finally, if $x$ is a complex number, we denote by $x^c$ the
complex conjugate of $x$.

\section{One dimensional mod \texorpdfstring{$pq$}{pq} representations.}\label{1-diml}

Given representations $\rho: G_K \rightarrow \overline\F_p^{\times}$ and $\rho': G_K \rightarrow \overline\F_q^{\times}$ we would
like to determine when $\rho$ and $\rho'$ arise from one Hecke character simultaneously with respect to the embeddings fixed
above. We know of no reasonable notion of finitely many mod $p$ representations with $p$ varying, such as $\rho$ and $\rho'$,
being compatible (see \cite{K-Pre}) but their arising from a Hecke character (see Section 4 of \cite{K-Pre}) is a special
circumstance which we will explore below.

\subsection{Hecke characters and Galois representations.}

In this subsection we wish to recall some properties of Hecke
characters, in particular properties of their attached mod $\ell$
Galois representations. We shall follow mainly the exposition of
\cite{S1}, section {\bf 3.4}.

\subsubsection{Hecke characters.}\label{A_0} A Hecke
character of $K$ is a continuous homomorphism $\chi : C_K
\rightarrow \C^{\times}$ whose restriction to the connected
component of 1 at infinity has form:$$x\mapsto \prod_{\sigma
\mbox{ real}} x_{\sigma}^{n_{\sigma}} \cdot \prod_{\sigma \mbox{
complex}} x_{\sigma}^{n_{\sigma}}
(x_{\sigma}^c)^{n_{\sigma^c}}$$with integers $n_{\sigma}$,
$\sigma \in \Sigma$. We have used $\sigma^c$ to denote the
element of $\Sigma$ given by $\sigma^c(x) = \sigma(x)^c$.
\medskip

We can consider the restriction of $\chi$ to $U_{\lambda}
\hookrightarrow C_K$ for finite primes $\lambda$. As $\chi$ is a
continuous homomorphism, almost all of these restrictions
$\chi|_{U_{\lambda}}$ are trivial and in any case have finite
orders and hence conductors. One defines the conductor ${\mathfrak
m}$ of $\chi$ to be the product of these local conductors.
\smallskip

As in \cite{S1}, we define for an arbitrary prime $\lambda$ of
$K$ the group $U_{\lambda,{\mathfrak m}}$: If ${\mathfrak m} =
\prod_{v\in S} \p_v^{m_v}$ with $S$ a finite set of finite primes
of $K$, the group $U_{\lambda,{\mathfrak m}}$ is defined to be
the connected component of $1\in K_{\lambda}^{\times}$ if
$\lambda$ is infinite, the whole of $U_{\lambda}$ if $\lambda
\not \in S$, and finally the group of units in $U_{\lambda}$ of
level $\ge m_{\lambda}$ if $\lambda \in S$. We may also define
$U_{\mathfrak m}$ as the product of all $U_{\lambda,{\mathfrak
m}}$.
\smallskip

The integers $n_{\sigma}$, $\sigma \in \Sigma$, are uniquely
determined by $\chi$. We can thus refer to the formal expression
$\sum_{\sigma \in \Sigma} n_{\sigma} \cdot \sigma \in \Z[\Sigma]$
as the $\infty$-type of $\chi$.
\smallskip

The classification of those elements in $\Z[\Sigma]$ that occur
as the $\infty$-type of {\it some} Hecke character is known, cf.
\cite{weil}, \cite{S}, but we shall not need to go into this.

\subsubsection{Galois representations attached to Hecke
characters.}\label{galois-hecke} Suppose that $\chi$ is a Hecke character of $\infty$-type $\sum_{\sigma \in \Sigma} n_{\sigma}
\cdot \sigma$ and conductor ${\mathfrak m}$. Let $\ell$ be a prime number. To $\chi$ is attached an $\ell$-adic representation
$\chi_{\ell}$ of $G_K$. The character $\chi_{\ell}$ takes its values in the group of units of $\overline\Q_{\ell}$, and as we
have fixed an isomorphism of the residue field of $\overline\Q_{\ell}$ with $\overline\F_{\ell}$, we can speak of the mod $\ell$
reduction $\chibar_{\ell}$ of $\chi_{\ell}$. We shall recall, cf. \cite{S1}, section {\bf 3.4}, the description of $\chi_{\ell}$
and $\chibar_{\ell}$ on inertia groups, i.e. on local units groups if we view these characters as characters on $C_K$.
\smallskip

Let $\lambda$ be a finite prime. Define the character
$\varepsilon_{\lambda}$ of $U_{\lambda}$ as the restriction
$\chi|_{U_{\lambda}}$. Then by construction of $\chi_{\ell}$, cf.
loc. cit., we have on $U_{\lambda}$:$$\chi_{\ell}(x) := \left\{
\begin{array}{ll} \varepsilon_{\lambda} (x) & , \mbox{
for } x\in U_{\lambda}, ~ \lambda \nmid \ell \\
\varepsilon_{\lambda} (x) \cdot \prod_{\sigma \in \Sigma(\lambda)}
\sigma_{\lambda} (x^{-1})^{-n_{\sigma}} & , \mbox{ for } x\in
U_{\lambda}, ~ \lambda \mid \ell . \end{array} \right.$$(Notice
that there are 2 different ways of establishing a correspondence
between the point of view of \cite{S1} and our view of a Hecke
character as a character of $C_K$. We could have changed signs on
the $n_{\sigma}$'s but would then have to define
$\varepsilon_{\lambda}$ as the {\it inverse} of the restriction of
$\chi$ to $U_{\lambda}$).

Notice, that as a consequence of definitions we have for a totally
positive global unit $u$ that$$\prod_{\lambda \mbox{ finite}}
\varepsilon_{\lambda}(u) = \prod_{\sigma \in \Sigma}
\sigma(u)^{-n_{\sigma}}.$$

The following lemma reverses this line of reasoning.

\begin{lemma}\label{hecke-lift} Suppose that for each finite prime
$\lambda$ of $K$ a continuous character$$\varepsilon_{\lambda}
\colon U_{\lambda} \longrightarrow \C^{\times}$$is given such
that only finitely many $\varepsilon_{\lambda}$ are non-trivial.
Suppose also that integers $n_{\sigma}$, $\sigma \in \Sigma$, are
given.
\smallskip

Then there exists a Hecke character $\chi$ on $C_K$ of
$\infty$-type $\sum_{\sigma \in \Sigma} n_{\sigma} \cdot \sigma$
whose restriction to $U_{\lambda}$ is $\varepsilon_{\lambda}$ for
all $\lambda$ if and only if$$\prod_{\lambda \mbox{ finite}}
\varepsilon_{\lambda} (u) = \prod_{\sigma \in \Sigma}
\sigma(u)^{-n_{\sigma}}$$for every totally positive global unit
$u$.
\end{lemma}
\begin{proof} As we saw above, the condition is necessary.

To prove sufficiency, consider the group$$U :=\prod_{\lambda
\mbox{ finite}} U_{\lambda}$$as a compact subgroup of $C_K$. The
$\varepsilon_{\lambda}$ define a continuous
character$$\varepsilon := \prod_{\lambda \mbox{ finite}}
\varepsilon_{\lambda}$$ on $U$. There is an extension $\chi$ of
$\varepsilon$ to a continuous character of $C_K$. The question is
whether there is an extension which has the desired shape on the
connected component $D_K$ of 1 in $C_K$.
\smallskip

Consider first the case where all $n_{\sigma}=0$. The question is then whether $\varepsilon$ extends to a continuous character of
$C_K$ vanishing on the closed subgroup $D_K$. This will be the case if and only if $\varepsilon$ vanishes on the intersection
$U\cap D_K$. As this intersection consists precisely of the embeddings into $U$ of totally positive global units, the desired
vanishing is implied by the condition we have imposed.
\smallskip

In the general case let ${\mathfrak m}$ be the conductor of
$\varepsilon$. A consequence of our hypothesis is then
that$$\prod_{\sigma \in \Sigma} \sigma(u)^{n_{\sigma}} =1$$for
every $u\in E_{\mathfrak m} := E\ \cap U_{\mathfrak m}$. This is
sufficient to ensure the existence of {\it some} Hecke character
$\chi$ on $C_K$ with $\infty$-type $\sum_{\sigma \in \Sigma}
n_{\sigma} \cdot \sigma$. The question is then whether the
character $\varepsilon \cdot (\chi|_U)^{-1}$ extends to a
character of finite order on $C_K$, and we are reduced to the
previous case.
\end{proof}

Returning now to the discussion before the lemma, we find for the
mod $\ell$ reduction $\chibar_{\ell}$ of the galois
representation $\chi_{\ell}$ that$$\chibar_{\ell} (x) =
(\varepsilon_{\lambda} ~\mathrm{mod} ~\ell)(x) \cdot \left\{
\begin{array}{ll} 1 & , \mbox{ for }
x\in U_{\lambda}, ~ \lambda \nmid \ell \\
\prod_{\sigma \in \Sigma(\lambda)} \sigmabar_{\lambda} ((x \mbox{
mod } \lambda)^{-1})^{-n_{\sigma}} & , \mbox{ for } x\in
U_{\lambda}, ~ \lambda \mid \ell . \end{array} \right.$$So, we
may view this as describing the galois representation
$\chibar_{\ell}$ on inertia. Viewing the homomorphisms$$(x \mbox{
mod } \lambda) \mapsto \sigmabar_{\lambda} ((x \mbox{ mod }
\lambda)^{-1})$$as characters on tame inertia above $\lambda$,
these characters are {\it fundamental characters} of level $n$ if
$\# k_{\lambda} = \ell^n$, with $\ell$ the residue characteristic
of $\lambda$, cf. again \cite{S1}. So, for each $\sigma \in
\Sigma(\lambda)$ we have a well-defined number
$\kappa(\lambda,\sigma) \in \{ 0, \ldots , n-1\}$ such
that:$$\sigmabar_{\lambda} =
\theta_{\lambda}^{\ell^{\kappa(\lambda,\sigma)}} .$$
\smallskip

Summing up, and viewing $\chibar_{\ell}$ as a character on $G_K$,
we have a description of $\chibar_{\ell}$ on inertia
groups:$$\chibar_{\ell} |_{I_{\lambda}} = \left\{
\begin{array}{ll} (\varepsilon_{\lambda} ~\mathrm{mod} ~\ell) & ,
\mbox{ for } \lambda \nmid \ell \\
(\varepsilon_{\lambda} ~\mathrm{mod} ~\ell) \cdot
\theta_{\lambda}^{-\sum_{\sigma \in \Sigma(\lambda)} n_{\sigma}
\cdot \ell^{\kappa(\lambda,\sigma)} } & , \mbox{ for } \lambda
\mid \ell \end{array} \right.$$with $\varepsilon_{\lambda}$ the
complex character obtained by restricting $\chi$ to
$U_{\lambda}$, viewed as a character on $I_{\lambda}$.

\subsection{Preliminary lemmas}

We begin by disposing off the easy case of {\it Artin lifts}.

\begin{lemma}\label{artin} Let $G$ be a profinite group and let $\tau: G
\rightarrow \overline\F_p^{\times}$ and $\tau': G \rightarrow
\overline\F_q^{\times}$ be two continuous (the target is given the
discrete topology), one-dimensional characters of $G$.
Then there is a complex character $\varepsilon :G \rightarrow
\C^{\times}$ which lifts $\tau$ and $\tau'$ (w.r.t. the
embeddings $\iota_p$ and $\iota_q$ fixed earlier) if and only if
for any two Artin lifts $\tilde\tau$ and $\tilde\tau'$ of $\tau$
and $\tau'$ the order of the character $\tilde\tau
\tilde{\tau'}^{-1}$ divides a power of $pq$.

If a lift $\varepsilon$ exists it is unique and in fact
$\varepsilon = \tilde\tau \psi' = \tilde\tau' \psi$ if
$\tilde\tau$ and $\tilde\tau'$ are Artin lifts and $\tilde\tau
\tilde{\tau'}^{-1} = \psi {\psi'}^{-1}$ with $\psi$ and $\psi'$
characters of $p$- and $q$-power order, respectively.
\end{lemma}
\begin{proof} This follows from the fact that the only roots of unity that
have trivial reduction mod $p$ are those that have order a power
of $p$.
\end{proof}

\begin{lemma}\label{artindescent}
Let $G$ be a profinite group and $G'$ a closed subgroup of finite index. Then for a (finite order, continuous) $\varepsilon': G'
\rightarrow \C^{\times}$ to be the restriction of a (finite order, continuous) character $\varepsilon: G \rightarrow \C^{\times}$
it is necessary and sufficient that $\ebar'_p$ and $\ebar'_q$, the reductions of $\varepsilon'$ mod $p$ and mod $q$, arise by
restriction from characters of $G$.
\end{lemma}
\begin{proof}
As $\ebar'_p$ arises by restriction of a mod $p$ character
$\ebar_p$ of $G$ it follows by choosing any lift of $\ebar_p$ to a
complex character that ${\varepsilon'}^{p^{\alpha}}$ arises by
restriction from $G$  for some non-negative integer $\alpha$.
Similarly ${\varepsilon'}^{q^{\beta}}$ arises by restriction from $G$.
Choosing integers $a$ and $b$ such that $ap^{\alpha}+bq^{\beta}=1$
we are done.
\end{proof}

\begin{rmk} It will be interesting to resolve the following question.
Let $\rho: G_K \rightarrow \overline\F_p^{\times}$ and $\rho': G_K \rightarrow \overline\F_q^{\times}$ be two continuous
one-dimensional characters of $G_K$. Let $L/K$ be a finite Galois extension of $K$ such that the restrictions $\rho|_{G_L}$ and
$\rho'|_{G_L}$ arise from a Hecke character $\chi'$ of $L$. Then do these restrictions also arise from $\chi \circ
\mathrm{Nm}_{L/K}$ for some Hecke character $\chi$ of $K$? The case when $\chi'$ is of finite order is dealt with by Lemma
\ref{artindescent}, and then $\chi'$ itself arises by base change from $K$ which certainly will not be true in general.
\end{rmk}

\begin{lemma}\label{necc} A necessary condition for $\rho$ and $\rho'$ to arise from
a Hecke character is that for all primes $\lambda$ of $K$ not above $p$ and $q$, and any two Artin lifts
$\tilde{\rho}|_{I_{\lambda}}$ and $\tilde{\rho'}|_{I_{\lambda}}$ of $\rho|_{I_{\lambda}}$ and $\rho'|_{I_{\lambda}}$ the order of
the character $\tilde{\rho}|_{I_{\lambda}} \tilde{\rho'}|_{I_{\lambda}}^{-1}$ divides a power of $pq$.
\end{lemma}
\begin{proof} This is clear from the construction of $p$-adic
representations from Hecke characters (see \cite{S}) and Lemma \ref{artin} (see Section 4 of \cite{K-Pre}). The point here is
that if $\chi$ is a Hecke character and if $\chi_p$ is the 1-dimensional $p$-adic representation of $G_K$ associated to it, then
the image of a inertia group at any prime $\lambda$ not above $p$ is finite and does not depend on the $p$ of the compatible
system of $p$-adic representations that $\chi$ gives rise to as long as the residue characteristics of $p$ and $\ell$ are
coprime.
\end{proof}

\begin{lemma}\label{inertia}
Assume the necessary conditions above. Then there is an Artin
character $\varepsilon:G_K \rightarrow \C^{\times}$ such that
$\ebar_p^{-1}\rho$ and $\ebar_q^{-1}\rho'$ are unramified outside
primes above $p$ and $q$ and the infinite places.
\end{lemma}
\begin{proof}
Let $L$ be the finite abelian extension of $K$ that is the
compositum of the fixed fields of the kernels of $\rho$ and
$\rho'$. Let $K'$ be the fixed field of $L$ under the action of
all the inertia groups $I_v$ for $v$ a finite place of $K$ not
above $p$ or $q$. As we are assuming the necessary conditions of
Lemma \ref{necc}, we deduce that there are lifts $\tilde{\rho}$
and $\tilde{\rho'}$ of $\rho|_{G_{K'}}$ and $\rho'|_{G_{K'}}$,
respectively, such that the order of $\tilde{\rho}
\tilde{\rho'}^{-1}$ is of the form $p^{\alpha}q^{\beta}$. Then by
Lemma \ref{artin}, there is an Artin character $\varepsilon'$ of
$G_{K'}$ such that $\ebar'_p=\rho|_{G_{K'}}$ and
$\ebar'_q=\rho|_{G_{K'}}$. Using Lemma \ref{artindescent}, we get
an Artin character $\varepsilon$ of $G_K$ such that
$\ebar_p=\rho|_{G_{K'}}$ and $\ebar_q=\rho|_{G_{K'}}$. Now
as $\ebar_p^{-1}\rho|_{G_{K'}}$ and $\ebar_q^{-1}\rho'|_{G_{K'}}$ are
trivial and $K'/K$ is unramified outside primes dividing $pq$ and the
infinite places we are done.
\end{proof}

\begin{lemma}\label{conductors}
If $\rho,\rho'$ as above are unramified outside primes above $p$
and $q$ and the infinite places, then any Hecke character $\chi$
that gives rise to both $\rho$ and $\rho'$ is of conductor that
divides $p^{\alpha}q^{\beta}$ for some non-negative integers
$\alpha,\beta$. More precisely, the conductor of any Hecke
character that gives rise to $\rho$ and $\rho'$ divides
${\mathfrak m}= \prod_{v_i} \mathrm{l.c.m} \{ \p_{v_i},
\mathrm{cond}(\rho'|_{I_{v_i}}) \} \prod_{v'_j} \mathrm{l.c.m} \{
\p_{v'_j}, \mathrm{cond}(\rho|_{I_{v'_j}}) \}$. Here, the $v_i$
and the $v'_j$ are the places of $K$ above $p$ and $q$
respectively.
\end{lemma}
\begin{proof}
This is again clear from the construction of $\lambda$-adic
representations from Hecke characters (see \cite{S} or the above
review in section \ref{galois-hecke}) and Lemma \ref{artin}, by
considering the restrictions $\rho|_{I_{v'_j}}$ and
$\rho'|_{I_{v_i}}$.
\end{proof}

Lemmas \ref{inertia} and \ref{conductors}
reduce us to examining $\rho,\rho'$ that are unramified
outside primes above $p$ and $q$
and the infinite places, and looking for lifts
by Hecke characters that are of conductor divisible only by primes
above $p$ and $q$.

\subsection{The case \texorpdfstring{$K=\Q$}{K=Q}}

Before studying the general case we first study the case $K=\Q$ that
will highlight the essential features.
Assume that $p,q$ are odd primes. As notation, for any prime $p$
we denote by $\theta_p := \theta_{\F_p}$ the mod $p$ cyclotomic
character and by $\tilde\theta_p$ its Teichm\"{u}ller lift (we
will regard the latter as a character taking values either in
$\overline\Q^{\times}$ or in $\overline\Q_p^{\times}$).

\begin{prop}\label{Q} Assume $K=\Q$, and assume the conditions of the
Lemma \ref{necc} above. Let $k_p$ be an integer defined mod $p-1$
such that $\rho|{I_p}=\theta_p^{k_p}$ and let $a_p$ be an integer
well-defined mod $A_p:=$ (the prime to $q$ part of $p-1$) such
that $\rho'|_{I_p}$ is the reduction mod $q$ of
$\tilde\theta_p^{a_p} \cdot \psi_p$ where $\psi_p$ is a character
of conductor and order a power of $p$. Define $k_q,a_p,B_q$
similarly (where this time for instance $B_q$ is the prime to $p$
part of $q-1$). The representations $\rho$ and $\rho'$ arise
simultaneously from a Hecke character if and only if there is an
integer $k$ that is simultaneously congruent to $k_p-a_p$ mod
$A_p$ and $k_q-a_p$ mod $B_q$. Note that the above conditions
depends only on
$\rho|_{I_p},\rho|_{I_q},\rho'|_{I_q},\rho'|_{I_p}$.
\end{prop}

\begin{proof}  We will assume
(w.l.o.g. because Lemma \ref{inertia} ensures this situation
after twisting by an Artin character)
above that $\rho$ and $\rho'$ are unramified outside $p,q$ and the
infinite place.

The Hecke characters for $\Q$ have a particularly simple
description: they are of the form $\varepsilon \cdot
\mathrm{Nm}^k$ where $\varepsilon$ is an Artin character,
$\mathrm{Nm}$ is the norm character, and $k$ is an integer. The
proof will be accomplished by analyzing mod $p$ and mod $q$
representations which arise from a character
$\chi:=\varepsilon\varepsilon'\mathrm{Nm}^k$ where $\varepsilon$
is primitive of conductor a power of $p$, and $\varepsilon'$ is
primitive of conductor a power of $q$ (we need only study such
characters because of Lemma \ref{conductors}). The mod $p$ and
mod $q$ (i.e., more precisely w.r.t. embeddings $\iota_p,\iota_q$
fixed above) characters of $G_K$ that arise from $\chi$ are
$\chibar_p:=\ebar_p \ebar'_p \theta_p^k$ and $\chibar_q:=\ebar_q
\ebar'_q \theta_q^k$.

We want to determine conditions on $\rho$ and $\rho'$ such that
there is a triple $(\varepsilon,\varepsilon',k)$ and the
corresponding Hecke character
$\chi=\varepsilon\varepsilon'\mathrm{Nm}^k$ has the property
$\chibar_p=\rho$ and $\chibar_q=\rho'$, i.e., gives rise to $\rho$
and $\rho'$. It will be useful first to determine to what extent
$\chibar_p$ and $\chibar_q$ determines the triple
$(\varepsilon,\varepsilon',k)$. For this note that
${\chibar_q}|_{I_p} ={\ebar_q}|_{I_p}$, and thus the Artin
character $\varepsilon$ which is of conductor a power of $p$ is
determined up to characters (of conductor a power of $p$) that
have order a power of $q$. This determines the wild part of
$\varepsilon$, and determines the tame part up to
characters of order a power of $q$.

By considering $\rho'|_{I_p}$ we get an integer $a_p$ that is
well-defined mod $A_p$ (the prime to $q$ part of $p-1$) such that
$\chibar_q|_{I_p}$ is the reduction mod the place above $q$
(fixed by $\iota_q$) of $\tilde\theta_p^{a_p}$ where now $a_p$
can vary mod $A_p$ (the reduction of a power of $\tilde\theta_p$
mod $q$ is constant when the power varies in a congruence class
mod $A_p$). A similar analysis applies to $\varepsilon'$ and
gives an integer $b_q$ that is well defined mod $B_q$ (the prime
to $p$ part of $q-1$) and such that $\chibar_p|_{I_q}$ is the
reduction mod the place above $p$ (fixed by $\iota_p$) of
$\tilde\theta_p^{b_q}$ where now $b_q$ can vary mod $B_q$.
Consider $\chibar_p|_{I_p}=\theta_p^{k_p}$ and
$\chibar_q|_{I_q}=\theta_q^{k_q}$ where $k_p$ and $k_q$ are
integers that are well-defined modulo $p-1$ and $q-1$
respectively by these equations.

We deduce from the above analysis, using further the fact that the
abelian extension $L$ of $\Q$ which is the compositum of the
fixed fields of the kernels of $\rho$ and $\rho'$ is generated by
the inertia groups above $p$ and $q$ in $\mathrm{Gal}(L/\Q)$
($\Q$ has no non-trivial unramified extensions), that there is a
Hecke character $\chi$  with $\chibar_p=\rho$ and
$\chibar_q=\rho'$, if and only if there exists an integer $k$ that
is congruent to $k_p-a_p$ mod $A_p$ and $k_q-b_q$ mod $B_q$. Thus
we have completely explicit necessary and sufficient conditions
for the existence of a $\chi$ such that $\chi_p=\rho$ and
$\chi_q=\rho'$ and these conditions depend only on
$\rho|_{I_p},\rho|_{I_q},\rho'|_{I_q},\rho'|_{I_p}$.
\end{proof}

\begin{rmk}
By examining the proof we see that we have a classification of
the lifts of $\rho$ and $\rho'$ by Hecke characters.
\end{rmk}

\subsection{The case of general \texorpdfstring{$K$}{K}}

We retain of course the notation of \ref{notation}, but also of
section \ref{galois-hecke} above.
\smallskip

The following theorem will give a general criterion for the
existence of a Hecke character $\chi$ of $K$, such
that$$\chibar_p \cong \rho \otimes \phi \qquad \mbox{and} \quad
\chibar_q \cong \rho' \otimes \phi'$$for some {\it unramified}
characters $\phi,\phi'$ on $G_K$.
\smallskip

Because of Lemma \ref{inertia} we may, and will, restrict
ourselves to the case where $\rho$ and $\rho'$ are both
unramified outside $p$ and $q$.
\medskip

Before stating the theorem we introduce some data attached to the
given representations $\rho$ and $\rho'$:
\smallskip

Let $v_1,\ldots ,v_s$ and $v'_1,\ldots ,v'_t$ be the primes of
$K$ above $p$ and $q$, respectively. Define the natural numbers
$A_i$, $B_j$ by the requirements:$$\# k_{v_i} -1 = A_i \cdot
(\mbox{power of } q), \quad q\nmid A_i,$$and similarly$$\#
k_{v'_j} -1 = B_j \cdot (\mbox{power of } p), \quad p\nmid B_j.$$

Let $k_i$ and $k'_j$ be integers such that:$$\rho|_{I_{v_i}} =
\theta_{v_i}^{k_i} , \quad i=1,\ldots ,s ,$$$$\rho'|_{I_{v'_j}} =
\theta_{v'_j}^{k'_j} , \quad j=1,\ldots ,t.$$

Also, there are integers $b_j$, well-defined modulo $B_j$, and a
complex character $\psi_j : I_{v'_j} \rightarrow \C^{\times}$ of
$q$-power order such that:$$(\rho|_{I_{v'_j}}) =
(\tilde\theta_{v'_j}^{b_j} ~\mathrm{mod} ~p)\cdot (\psi_j
~\mathrm{mod} ~p) ,$$and similarly we have integers $a_i$,
well-defined modulo $A_i$, such that:$$(\rho'|_{I_{v_i}}) =
(\tilde\theta_{v_i}^{a_i} ~\mathrm{mod} ~q) \cdot (\psi'_i
~\mathrm{mod} ~q),$$with a complex character $\psi'_i : I_{v_i}
\rightarrow \C^{\times}$ of $p$-power order.
\medskip

Finally, we shall denote by $\mathrm{rec}_{\lambda}$ the
reciprocity map $K_{\lambda}^{\times} \rightarrow
\mathrm{Gal}(K_{\lambda}^{ab}/K_{\lambda})$.

\begin{theorem}\label{criterion} Consider the above situation and
suppose that integers $n_{\sigma}$, $\sigma \in \Sigma$, are
given. Then there exists a Hecke character $\chi$ of
$\infty$-type $\sum_{\sigma \in \Sigma} n_{\sigma} \cdot \sigma$,
and unramified characters$$\phi \colon G_K \longrightarrow
\overline\F_p^{\times}, \qquad \phi' \colon G_K \longrightarrow
\overline\F_q^{\times},$$such that$$\chibar_p \cong \rho \otimes
\phi \qquad \mbox{and} \quad \chibar_q \cong \rho' \otimes \phi'
,\leqno({\ast})$$if and only if the following conditions hold:

$$k_i - a_i \equiv -\sum_{\sigma \in
\Sigma(v_i)} n_{\sigma} \cdot p^{\kappa(v_i,\sigma)} =: \xi_i
\quad \mathrm{mod} ~~A_i,\qquad \mbox{for $i=1,\ldots ,s$,}
\leqno{(1)}$$and

$$k'_j - b_j \equiv -\sum_{\sigma \in
\Sigma(v_j)} n_{\sigma} \cdot q^{\kappa(v_j,\sigma)} =:\xi'_j
\quad \mathrm{mod} ~~B_j,\qquad \mbox{for $j=1,\ldots ,t$.}
\leqno{(1')}$$

\noindent $(2)$ We have$$\prod_i (\tilde\theta_{v_i}^{k_i -
\xi_i} \cdot \psi'_i) \circ \mathrm{rec}_{v_i}(u) \cdot \prod_j
(\tilde\theta_{v'_j}^{k'_j - \xi'_j} \cdot \psi_j) \circ
\mathrm{rec}_{v'_j}(u) = \prod_{\sigma \in \Sigma}
\sigma(u)^{-n_{\sigma}},$$for every totally positive global unit
$u$.
\end{theorem}
\begin{proof} The structure of the proof is as follows. We first
seek to determine the restriction of $\chi$ to local unit groups
$U_{\lambda}$ for primes $\lambda$ lying above $p$ or $q$. This
is done via matching the description we have from section
\ref{galois-hecke} with the information coming from the behaviour
of $\rho$ and $\rho'$ at local inertia groups. Then Lemma
\ref{artin} is invoked to give local conditions which will turn
out as the conditions $(1)$ and $(1')$. Once these restrictions
of $\chi$ have been determined, the existence of $\chi$ amounts
to the global condition of Lemma \ref{hecke-lift} which will turn
into condition $(2)$.
\medskip

\noindent {\it Proof of necessity:} Assume that $\chi$ exists with
the stated properties. According to lemma \ref{conductors},
$\chi$ is unramified outside the primes dividing $p$ or $q$, and
the infinite primes. If $\ell$ is a prime number and $\lambda$ a
finite prime of $K$, we have -- retaining the notation of section
\ref{galois-hecke} -- that$$\chibar_{\ell} |_{I_{\lambda}} =
\left\{ \begin{array}{ll} (\varepsilon_{\lambda} ~\mathrm{mod}
~\ell) & , \mbox{ for } \lambda \nmid \ell \\
(\varepsilon_{\lambda} ~\mathrm{mod} ~\ell) \cdot
\theta_{\lambda}^{-\sum_{\sigma \in \Sigma(\lambda)} n_{\sigma}
\cdot \ell^{\kappa(\lambda,\sigma)} } & , \mbox{ for } \lambda
\mid \ell \end{array} \right.$$ In particular we have for
$i=1,\ldots ,s$:$$\chibar_q |_{I_{v_i}} = (\varepsilon_{v_i}
~\mathrm{mod} ~q)= \rho' |_{I_{v_i}} = (\tilde\theta_{v_i}^{a_i}
~\mathrm{mod} ~q) \cdot (\psi'_i ~\mathrm{mod}
~q)$$and$$\chibar_p |_{I_{v_i}} = (\varepsilon_{v_i}
~\mathrm{mod} ~p) \cdot \theta_{v_i}^{\xi_i} = \rho |_{I_{v_i}} =
\theta_{v_i}^{k_i}.$$
\smallskip

We conclude that the mod $p$ and mod $q$
characters$$\theta_{v_i}^{k_i-\xi_i} \qquad \mbox{and} \qquad
(\tilde\theta_{v_i}^{a_i} ~\mathrm{mod} ~q) \cdot (\psi'_i
~\mathrm{mod} ~q)$$of $I_{v_i}$ simultaneously lift to a
(uniquely determined) complex character of $I_{v_i}$, namely
$\varepsilon_{v_i}$. On the other hand, as individual lifts of
these characters are$$\tilde\theta_{v_i}^{k_i-\xi_i} \qquad
\mbox{and} \qquad \tilde\theta_{v_i}^{a_i} \cdot \psi'_i$$where
$\psi'_i$ has $p$-power order, and as $\tilde\theta_{v_i}$ has
order prime to $p$, we deduce from lemma \ref{artin} that the
character$$\tilde\theta_{v_i}^{a_i - k_i + \xi_i}$$ has $q$-power
order. This implies $(1)$ by the definition of $A_i$. Lemma
\ref{artin} also gives:$$\varepsilon_{v_i} =
\tilde\theta_{v_i}^{k_i - \xi_i} \cdot \psi'_i. \leqno{(\ast
\ast)}$$
\smallskip

Similarly, $(1')$ follows from the fact $\chi$ gives a
simultaneous lift of the representations $\rho$ and $\rho'$
restricted to $I_{v'_j}$, and we see that in
fact$$\varepsilon_{v'_j} = \tilde\theta_{v'_j}^{k'_j - \xi'_j}
\cdot \psi_j. \leqno{(\ast \ast \ast)}$$
\smallskip

Now, as $\chi$ is unramified outside $pq$ and the infinite
primes, condition $(2)$ follows from lemma \ref{hecke-lift},
$(\ast \ast)$ and $(\ast \ast \ast)$.
\medskip

\noindent {\it Proof of sufficiency:} Define complex characters on
inertia at the primes $v_i$, $i=1,\ldots ,s$, and $v'_j$,
$j=1,\ldots ,t$, according to $(\ast \ast)$ and $(\ast \ast
\ast)$ above. Define also $\varepsilon_{\lambda}$ to be the
trivial character on $I_{\lambda}$ if $\lambda$ is a finite prime
not above $p$ or $q$.

If these characters are viewed as characters on local unit groups
$U_{\lambda}$, then $(2)$ combined with Lemma \ref{hecke-lift}
imply that these local characters are in fact restrictions to the
$U_{\lambda}$ of a Hecke character $\chi$ with $\infty$-type
$\sum_{\sigma \in \Sigma} n_{\sigma} \cdot \sigma$.
\smallskip

For such a $\chi$ we deduce, utilizing $(1)$ and $(1')$, and
reversing the pertinent reasoning in the proof of necessity
above, that the characters $\chibar_p^{-1} \cdot \rho$ and
$\chibar_q^{-1} \cdot \rho'$ vanish on every inertia group, and
hence are both globally unramified.
\end{proof}

\begin{rmk} As will be seen from the proof of the theorem, the
contribution at the place $v_i$ to the conductor of a lift $\chi$
is `essentially'$$\mathrm{l.c.m} \{ \p_{v_i}, \mathrm{cond}
(\psi'_i) \} = \mathrm{l.c.m} \{ \p_{v_i}, \mathrm{cond}
(\rho'|_{I_{v_i}}) \}.$$ We say `essentially' because this is
definitely true if $\psi'_i$ is nontrivial. If however $\psi'_i$
is trivial the conductor is 1 or $\p_{v_i}$ depending on whether
$k_i-\xi_i$ is divisible by $\# k_{v_i} -1$ or not; however, this
condition depends on the $n_{\sigma}$'s. Similar remarks apply to
the places dividing $q$ of course.
\end{rmk}

\begin{rmk} We chose to formulate a theorem considering only
representations unramified outside $pq$ and the infinite primes.
The reason is first of all that this is the crucial case, and
secondly that we wanted to avoid an unnecessarily complicated
statement. However, the observant reader will notice that our
reduction to the above case via lemma \ref{inertia} is not
completely explicit: The behavior of an $\varepsilon$ from lemma
\ref{inertia} at primes above $p$ or $q$ is not given explicitly.
\smallskip

However, this situation can easily be rectified, and we shall
limit ourselves to some indications: If one wants to prove a
statement as in theorem \ref{criterion} for arbitrarily ramified
representations, there are 2 possible, essentially equivalent,
ways of accomplishing this. First, an inspection of the proof of
theorem \ref{criterion} reveals immediately that the method
extends without any difficulty to the general case of arbitrary
ramification: First, one would have conditions like $(1)$, $(1')$
for each ramification point $w$ not dividing $pq\cdot \infty$;
these are purely local conditions coming directly from lemma
\ref{artin}. Secondly, condition $(2)$ would change by
incorporating certain factors on the left hand side; these
factors would depend on the restrictions $\rho|_{I_w}$ and
$\rho'|_{I_w}$ for each $w$ as above. We shall leave it to the
reader to work out the exact form of these factors.

The second way to extend theorem \ref{criterion} to the case of
arbitrary ramification is to be more precise about the behavior
of an $\varepsilon$ as in lemma \ref{inertia} at the primes
dividing $pq$. A possible way to do this would be to use lemma
\ref{hecke-lift} for the case $n_{\sigma}=0$ -- under the
assumptions of lemma \ref{inertia} -- to `shift' (via twisting)
any ramification outside $pq\cdot \infty$ to a prime above $p$
(say). After this, one could plug the data into theorem
\ref{criterion} and get explicit conditions. As the reader will
however quickly ascertain, this method would result in exactly
the same modified conditions $(1)$, $(1')$ and $(2)$ that result
from the first method.
\end{rmk}

\begin{rmk} \label{solutions} Notice that if $\chi_1$ and $\chi_2$
are 2 Hecke characters with the same $\infty$-type which are both
`solutions' to $(\ast)$, i.e. have their mod $p$ and mod $q$
representations isomorphic up to twist by unramified to $\rho$
and $\rho'$, respectively, then $\chi_1$ and $\chi_2$ differ by a
globally unramified Artin character. For the Hecke character
$\varepsilon := \chi_1 \chi_2^{-1}$ has trivial $\infty$-type,
i.e. is an Artin character, and its mod $p$ and mod $q$
reductions are by hypothesis both unramified whence the claim by
Lemma \ref{artin}.
\end{rmk}

\begin{rmk} The reason why we consider only `lifting up to
unramified' as in the Theorem above is to be found in the nature
of the conditions we impose: As an analysis of the argument
quickly reveals, in order to discuss {\it exact} lifting, and not
just `up to unramified', we would have to consider some sort of
analogue to Lemma \ref{hecke-lift} but with given restrictions to
decomposition groups instead of inertia groups. This ties up with
the Grunwald-Wang theorem where of course the problem is that one
cannot in general lift without introducing additional
ramification. We have not been able to find a reasonable
characterization of pairs $(\rho,\rho')$ with `exact' lifts and
doubt whether there is one. The following example -- a simple
counting argument -- shows that in general one cannot lift an arbitrary
pair of unramified characters.
\end{rmk}

\noindent {\bf Example:} Suppose that $K$ is an imaginary
quadratic field containing no other roots of unity that $\pm 1$.
Let $p$ and $q$ be odd primes that split in $K$. Fix an embedding
$\sigma : K\hookrightarrow \C$ and consider an $\infty$-type
$m\cdot \sigma + n\cdot \sigma^c$, with integers $m$ and $n$. We
apply Theorem \ref{criterion} to find the conditions on $m$ and
$n$ necessary and sufficient for the existence of a Hecke
character $\chi$ with $\infty$-type $m\cdot \sigma + n\cdot
\sigma^c$ and both $\chibar_p$ and $\chibar_q$ globally
unramified. So, we apply Theorem \ref{criterion} for the case
$\rho = \rho' =1$ and hence $k_i=k'_j=a_i=b_j=0$, and $\psi_i =
\psi'_j = 1$ for all $i,j$.
\smallskip

The conditions $(1)$ and $(1')$ of the theorem amount then to 4
congruence conditions on $n$ and $m$ which are seen to boil down
to:$$n,m \equiv 0 ~\mathrm{mod}~ C$$where $C:=\mathrm{l.c.m.} \{
A ,B\}$ with $A:=A_1=A_2$ and $B:=B_1=B_2$ as in the definitions
preceding the theorem. Now, $C$ is an even number so if $n$ and
$m$ are both divisible by $C$, then condition $(2)$ of the
theorem is automatically satisfied as $\pm 1$ are the only units
in $K$.

So there exist Hecke characters $\chi_1$ and $\chi_2$ with
$\infty$-types $C\cdot \sigma$ and $C\cdot \sigma^c$,
respectively, such
that$$((\overline{\chi_i})_p,(\overline{\chi_i})_q) =
(\phi_i,\phi'_i), \qquad i=1,2$$where $\phi_i, \phi'_i$, $i=1,2$,
are unramified characters. Now, if $\chi$ is any Hecke character
such that $\chibar_p$ and $\chibar_q$ are both unramified, its
$\infty$-type has shape $\mu C\cdot \sigma + \nu C\cdot \sigma^c$
with integers $\mu$ and $\nu$. Consequently, $\chi = \varepsilon
\chi_1^{\mu} \chi_2^{\nu}$ with $\varepsilon$ an Artin character
which must be globally unramified (cf. Remark \ref{solutions}
above).
\smallskip

So any pair $(\rho,\rho')$ of unramified mod $p$ and mod $q$ characters that lifts {\it exactly} to a Hecke character necessarily
has form:$$(\rho,\rho') = (\ebar_p \phi_1^{\mu} \phi_2^{\nu} ,\ebar_q (\phi'_1)^{\mu} (\phi'_2)^{\nu})$$with an unramified Artin
character $\varepsilon$. If the class group of $K$ has exponent $\alpha$ and order $h$, it follows that the number of such pairs
$(\rho,\rho')$ is bounded by $\alpha^2\cdot h$. On the other hand, if $h$ is not divisible by $p$ or $q$, then the total number
of pairs of unramified mod $p$ and mod $q$ characters is at least $h^2$. Hence, if further $h>\alpha^2$ there necessarily exist
at least 1 such pair that does not lift `exactly' to a Hecke character. As a concrete example, one can take $K=\Q(\sqrt{-3\cdot
5\cdot 7\cdot 11})$, which has class group isomorphic to $(\Z/2\Z )^3$, and for $p$ and $q$ any two distinct odd primes that
split in $K$.

\section{Serre's conjecture mod \texorpdfstring{$pq$}{pq}}\label{serre}

This will be a speculative section on issues surrounding
mod $pq$ versions of Serre's conjectures in \cite{Serre-Duke}.
Let $p$ and $q$ be odd primes and fix as above embeddings $\iota_p$
and $\iota_q$ of  $\overline\Q$ in $\overline\Q_p$ and
$\overline\Q_q$ as before. Barry Mazur and Ken Ribet had
investigated the question of when a pair of odd continuous
irreducible mod $p$ and mod $q$ Galois representations
$$\rho:G_{\mathbb Q} \rightarrow \GL_2(\overline{\mathbb F}_p)$$ and
$$\rho':G_{\mathbb Q} \rightarrow \GL_2(\overline{\mathbb F}_q)$$ arise from a
newform $f \in S_2(\Gamma_1(N))$ for some level $N$ and with
respect to the embeddings $\iota_p$ and $\iota_q$ that have been
fixed. William Stein made some computations (\cite{Stein}) towards this question
in which the level $N$ was varied in order to discover a newform
$f \in S_2(\Gamma_0(N))$ of the desired sort.

There are certain local constraints that one must impose, and
which we will come to later (see Conjecture \ref{serrepq} below),
to expect to have an affirmative answer. It must be pointed out
here that even assuming Serre's original conjecture
\cite{Serre-Duke}, i.e., assuming that $\rho$ and $\rho'$ are
individually modular, proving that they are {\it simultaneously}
modular is a different ball game. Even after having imposed these
necessary local constraints what makes this question difficult to
address even computationally is that there seems no a priori way
of guessing what levels $N$ one should be looking at  (though one
does know by the analysis in \cite{Carayol-Duke} what levels one
should not be looking at as local constraints preclude some
primes diving the levels of newforms which give rise to either of
the $\rho$ or $\rho'$!). Unlike the situation in Serre's
conjecture where he did specify a minimal level at which a
representation like $\rho$ should be found (before Ribet et al,
cf. \cite{Ribet-Mot}, proved this $\epsilon$-conjecture one did
not even know whether there existed a minimal level even without
caring about what it exactly was), there does not seem to be a
``minimal level'' $N_{\rho,\rho'}$ at which $\rho$ and $\rho'$
should arise from a newform $f \in
S_2(\Gamma_1(N_{\rho,\rho'}))$, i.e., all other levels $M$ at
which $\rho$ and $\rho'$ arise from a newform in
$S_2(\Gamma_1(M))$ should be divisible by $N_{\rho,\rho'}$.

There is one natural guess for which levels $N$ one can look at
which unfortunately is very unlikely to be correct. Let
$N(\rho)$, resp., $N(\rho')$, be the prime to $p$, resp., prime
to $q$, part of the Artin conductor of $\rho$, resp. $\rho'$.
Then we can try to look for the desired newform from which
$\rho,\rho'$ arise simultaneously in
$S_2(\Gamma_1(N(\rho)N(\rho')p^{\alpha}q^{\beta}))$ where
$\alpha,\beta$ vary. The reason why this is unlikely to work is
that there are only {\bf finitely many} $\alpha,\beta$ such that there
will exist a newform $f \in
S_2(\Gamma_1(N(\rho)N(\rho')p^{\alpha}q^{\beta}))^{new}$ that can
give rise both $\rho$ and $\rho'$. This follows from the analysis
in \cite{Carayol-Duke} which gives that if $\rho$ (resp.,
$\rho'$) arises from $S_2(\Gamma_1(N(\rho)N(\rho')
p^{\alpha}q^{\beta}))^{new}$ then $\beta$ (resp., $\alpha$) is
bounded. (Note that on the contrary by Lemme 1 of
\cite{Carayol-Duke} $\rho$, resp. $\rho'$, does arise from
$S_2(\Gamma_1(N(\rho)p^{\alpha}))^{new}$, resp.,
$S_2(\Gamma_1(N(\rho')q^{\beta}))^{new}$ for all $\alpha \geq 2$,
resp., $\beta \geq 2$.) Further again unlike the situation for a
single representation $\rho$, if there does exist a newform $f$
which gives rise to $\rho$ and $\rho'$, it is not clear that
there are infinitely many, i.e., there is no known way that one
can systematically raise levels mod $pq$ unlike in the case of mod
$p$ representations where a complete study is available because
of the work of Ribet (cf. \cite{Ribet-ICM}) and others.

In the questions raised by Mazur and Ribet and the subsequent
computations of William Stein only the level aspect of this
question was considered. In the following paragraphs we discuss
the ``weight aspect''.

Suppose for simplicity for this paragraph that $\rho$ and $\rho'$ arise individually
(with respect to the embeddings $\iota_p,\iota_q$)
from newforms of level $N$ in $S_k(\Gamma_0(N))$ with $N$
squarefree and $(N,pq)=1$. (This
ensures that the necessary local conditions of Conjecture \ref{serrepq}
below are satisfied.)
Then we can ask if they arise {\it simultaneously} from a newform $f \in
S_{k'}(\Gamma_0(N))$ for some $k'>>0$. Unlike the previous
paragraph this is more difficult to rule out as for all $k'$
congruent to $k$ mod $\mathrm{l.c.m.}(p-1,q-1)$ the $\rho$ and
$\rho'$ {\it individually} do arise from $S_{k'}(\Gamma_0(N))$.
This of course can be seen by multiplication by powers of the
normalized Eisenstein series $E_{p-1}$ and $E_{q-1}$ (Hasse
invariants mod $p$ and mod $q$) that are congruent to 1 mod $p$
and 1 mod $q$ respectively and then applying the mod $p$ or mod
$q$ Deligne-Serre lifting lemma, or better still in the present context
by multiplication by a suitable power of the souped up ``mod $pq$
Hasse invariant''
$$E_{\mathrm{l.c.m.}(p-1,q-1)}=1-{ \frac{ 2(\mathrm{l.c.m.}(p-1,q-1)) }
{ B_{\mathrm{l.c.m.}(p-1,q-1)} } }.\sum_{n \geq 1}
\sigma_{\mathrm{l.c.m.}(p-1,q-1)-1}(n)q^n$$ (where by $B_k$ we
mean the $k$th Bernoulli number and $\sigma_k(n)=\sum_{d \mid
n}d^k$) which by the von-Staudt-Clausen theorem is congruent to 1
mod $pq$. Unfortunately there is no mod $pq$ Deligne-Serre
lifting lemma, and so we do not know and perhaps do not even
expect that there is any periodicity with respect to weight of
newforms in $S_k(\Gamma_0(N))$ which give rise to $\rho$ and
$\rho'$ simultaneously.

It will be interesting to investigate the
issues raised in the paragraph computationally. What makes the
possibility of $\rho$ and $\rho'$ arising simultaneously from
$S_{k'}(\Gamma_0(N))$ attractive is that arithmetically the
weight is just one parameter (or automorphically, the component at
the archimedean place) while if we were to look at
$S_k(\Gamma_0(N))$ for fixed $k$ and varying $N$ we would {\it
have} to allow $N$ to be divisible by varying primes (as follows
from the analysis in \cite{Carayol-Duke}), so arithmetically vary
infinitely many independent parameters (or automorphically, vary
infinitely many non-archimedean local components).

In \cite{K-MRL} one set of local constraints for $\rho$ and
$\rho'$ to arise simultaneously from a newform were highlighted:
this may be called the {\it semistable case}.  We will
formulate here a finer mod $pq$ Serre conjecture.

Given a (algebraic) Weil-Deligne parameter $(\tau_{\ell},N_{\ell})$ (for this see \cite{Tate-Corvallis}) with
$\tau_{\ell}:W_{\Q_{\ell}} \rightarrow \GL_2(\overline\Q)$ a continuous representation of $W_{\Q_{\ell}}$, the Weil group of
$\Q_{\ell}$, with the target given the discrete topology, and $N_{\ell}$ a nilpotent matrix in $M_2(\overline{\Q})$, and $p\neq
\ell$ is a prime then the corresponding 2-dimensional $p$-adic representation $(\tau_{\ell},N_{\ell})_p:{\rm
Gal}(\overline\Q_{\ell}/\Q_\ell) \rightarrow \GL_2(\overline\Q_p)$ is:

\begin{itemize}
\item $N_{\ell}=0$, $\tau_{\ell}$ reducible:
$\tau_{\ell}=\varepsilon_1 \oplus \varepsilon_2$, the sum of a pair of
quasicharacters of the Weil group $W_{\Q_{\ell}}$ taking values in
$\overline\Q^{\times}$, then $(\tau_{\ell},N_{\ell})_p:=\varepsilon_{1,p}
\oplus \varepsilon_{2,p}$, and
$\varepsilon_{1,p},\varepsilon_{2,p}$ the corresponding $p$-adic
characters arising via the fixed embedding $\iota_p:\overline\Q
\rightarrow \overline\Q_p$.

\item $N_{\ell}=0$, $\tau_{\ell}$ irreducible: then up to twisting by a character
$\varepsilon:W_{\Q_{\ell}} \rightarrow \overline\Q^{\times}$,
$\tau_{\ell}$ has finite image. Then $(\tau_{\ell},N_{\ell})_p:= (\tau_{\ell} \otimes
\varepsilon)_p \otimes \varepsilon_p^{-1}$ with $(\tau_{\ell} \otimes
\varepsilon)_p$ the embedding of the finite image of $\tau_{\ell} \otimes
\varepsilon$ in $\GL_2(\overline\Q)$ into
$\GL_2(\overline{\Q_p})$ and  $\varepsilon_p$ the
$p$-adic character corresponding to $\varepsilon$ arising via the
fixed embedding $\iota_p:\overline\Q \rightarrow \overline\Q_p$.

\item $N_{\ell} \neq 0$: then $\tau_{\ell}$ is forced to be of the form $\varepsilon \oplus
\varepsilon || \ \ ||$ with $\varepsilon$ a quasicharacter of the
Weil group $W_{\Q_{\ell}}$ taking values in $\overline\Q^{\times}$ and
$|| \ \ ||$ the ``norm character'' (i.e., unramified character
with $||{\rm Frob}_{\ell}||=\ell$) then $(\tau_{\ell},N_{\ell})_p$
is $\left(\begin{matrix} \varepsilon_p||\
\ ||_p & 0 \\  0 & \varepsilon_p  \end{matrix}\right) {\rm exp}(N_{\ell}t_{p})$ :
here $\varepsilon_{p}$ is the $p$-adic character corresponding to
$\varepsilon$ w.r.t. the embedding
$\iota_p:\overline\Q \rightarrow \overline\Q_p$ as before,
$t_{p}$ is the homomorphism of the Galois group of the maximal
tamely
ramified extension of ${\bf Q}_{\ell}$ which arises by projecting
first to $I_{\ell}^t$, the
tame inertia at $\ell$ (i.e., after choosing a Frobenius $\sigma$,
$\sigma^n\tau \rightarrow \tau$ where $\tau \in I_{\ell}^t$),
followed by a $\Z_p$-valued homomorphism of
$I_{\ell}^t$ which is the projection to the unique
$\Z_p$ quotient of $I_{\ell}^t$ and $|| \ \ ||_p$ is the $p$-adic
cyclotomic character (the isomorphism class of $(\tau_{\ell},N_{\ell})_p$ is
independent of the choices made).
\end{itemize}

\begin{conjecture}\label{serrepq}
Let $\rho:G_{\mathbb Q} \rightarrow \GL_2(\overline{\mathbb
F}_p)$ and $\rho':G_{\mathbb Q} \rightarrow
\GL_2(\overline{\mathbb F}_q)$ be continuous, odd, absolutely
irreducible representations and assume that $\rho$ is unramified
at $q$ and $\rho'$ is unramified at $p$. Let $k(\rho)$ and
$k(\rho')$ be the Serre invariants of $\rho,\rho'$ as in
\cite{Serre-Duke}. Then $\rho$ and $\rho'$ are the reductions mod
$p$ and mod $q$ of the $p$-adic and $q$-adic representations
(w.r.t. $\iota_p,\iota_q$) attached to a newform in
$S_k(\Gamma_1(N))$ for some integer $k \geq 2$ and some positive
integer $N$ with $(N,pq)=1$ if and only if \begin{itemize} \item
there is an integer congruent to $k(\rho)$ mod $p-1$ and
$k(\rho')$ mod $q-1$, \item and further for every prime $\ell
\neq p,q$, there is a (algebraic) Weil-Deligne parameter
$(\tau_{\ell},N_{\ell})$ and a choice of integral models for $(\tau_{\ell},N_{\ell})_p$ and
$(\tau_{\ell},N_{\ell})_q$ such that their mod $p$ and mod $q$ reductions are
isomorphic to $\rho|_{D_{\ell}}$ and $\rho'|_{D_{\ell}}$
respectively.\end{itemize}
\end{conjecture}

\noindent{\bf Discussion of the conjecture:} The main feature of
the conjecture is that the conditions
are local (recall that in \cite{Serre-Duke} the definition of
$k(\rho),k(\rho')$ is purely in terms of $\rho|_{I_p},\rho'|_{I_q}$).
The condition in the conjecture at primes $\ell$
is empty for $\ell$ that are $\neq p,q$
and unramified in both $\rho$ and $\rho'$, and thus there
are only finitely many conditions that need to be checked
as per the conjecture. The restriction that $\rho$ is
unramified at $q$ and $\rho'$ is unramified at $p$ is a
technical restriction which arises from the difficulties
of studying mod $p$ representation restricted to $D_p$ arising
from newforms with levels divisible by powers of $p$ (there are some
results towards this in \cite{BM} and \cite{Ulmer}):
it will be nice to remove this technical restriction
and refine the above conjecture to allow $\rho$ and $\rho'$ to be
ramified at $q$ and $p$ respectively.
\smallskip

As the referee has remarked, in general the local constraints in the conjecture are non-empty: We can find distinct primes
$\ell$, $p$ and $q$, and representations $\rho \colon G_{\Q} \rightarrow \GL_2(\overline\F_p)$, $\rho' \colon G_{\Q} \rightarrow
\GL_2(\overline\F_q)$ such that $\rho$ has conductor divisible exactly by $\ell^3$ whereas $\rho'$ is unramified at $\ell$. The
results due to Carayol \cite{Carayol-Duke} show that the pair $(\rho,\rho')$ definitely can not arise from a modular form.
\smallskip

We could make a more optimistic conjecture and {\it fix} a weight $k \geq 2$ congruent to $k(\rho)$ mod $p-1$ and $k(\rho')$ mod
$q-1$, and ask that $\rho$ and $\rho'$ arise from $S_k(\Gamma_1(N))$ for some (variable) $N$ prime to $pq$ (as in \cite{Stein}),
or we can {\it fix} $N$ prime to $pq$ and divisible by $N(\rho)$ and $N(\rho')$, the prime-to-$p$ and prime-to-$q$ parts of the
Artin conductors of $\rho$ and $\rho'$ respectively, and ask that $\rho$ and $\rho'$ arise from $S_k(\Gamma_1(N))$ for some
(variable) $k$ as in the discussion above. We can also weaken the conjecture and instead of asking that $\rho$ and $\rho'$ arise
from a common newform with respect to fixed embeddings $\iota_p,\iota_q$, allow these embeddings to vary, or alternatively ask
that $\rho$ and $\rho'$ arise simultaneously (w.r.t. the fixed embeddings $\iota_p$ and $\iota_q$) from the {\it Galois orbit} of
a newform $f$.

One can also formulate a conjecture entirely analogous to Conjecture
\ref{serrepq} for compatible lifts of $\rho$ and $\rho'$ to
$p$-adic and $q$-adic representations of $G_{\Q}$, which we leave
as an exercise for the interested reader. As was pointed out in
\cite{K-MRL} such a conjecture has relevance to Serre's {\it
original} conjecture in \cite{Serre-Duke} because of the
modularity lifting theorems of Wiles et al.
\smallskip

In connection with Conjecture \ref{serrepq} it is important to
notice that the weight $k$ may very well depend on the choice of
primes of $\overline\Q$ over $p$ and $q$. The following simple
example illustrates this.
\smallskip

\noindent {\bf Example:} Suppose that $f$ and $f'$ are two
newforms on $\mathrm{SL}_2(\Z)$ of some weight $k$ and let $\rho$
and $\rho'$ be the mod $p$ and mod $q$ representations attached
to $f$ and $f'$, respectively. The conditions of Conjecture
\ref{serrepq} are then satisfied so that the Conjecture would
have $\rho$ and $\rho'$ arising from a newform $F\in
S_m(\Gamma_1(N))$ for some $m\ge 2$ and some $N$ prime to $pq$.

To illustrate with a simple example that the possible weights $m$
may depend on our choice of primes over $p$ and $q$, suppose that
$p=5$, $q=7$, and that $f$ and $f'$ are the 2 newforms of weight
24 and level 1:$$f = 24\alpha \Delta^2 + Q^3 \Delta \qquad
\mbox{and} \qquad f' = 24\alpha' \Delta^2 + Q^3 \Delta$$with $\{
\alpha, \alpha' \} := \{ -13/2 \pm \sqrt{144169}/2 \}$. The
primes 5 and 7 both split in the quadratic field
$\Q(\sqrt{144169})$, say $5={\mathfrak p}_5\cdot {\mathfrak
p}_5'$ and $7={\mathfrak p}_7\cdot {\mathfrak p}_7'$. Now,
$\Delta$ is congruent mod ${\mathfrak p}_7$ to either $f$ or
$f'$. Let us assume notation so that $\Delta \equiv f ~({\mathfrak
p}_7)$. As $Q\equiv 1 ~(5)$ we have then$$\Delta \equiv f
~({\mathfrak p}_5), \qquad \Delta \equiv f' ~({\mathfrak
p}_7').$$On the other hand, as $\alpha \alpha'$ is divisible by 5
we have$$f' \equiv f ~({\mathfrak p}_5), \qquad f' \equiv f'
~({\mathfrak p}_7)$$and $f'$ is the solution to these congruences
with smallest possible weight.

\begin{rmk}
Unlike the case of 1-dimensional mod $pq$ representations
(see Lemma \ref{conductors}) one would expect that there
are modular forms of levels divisible by arbitrarily many primes
that give rise to $\rho$ and $\rho'$ if the local conditions
are met (this is substantiated by the tables in \cite{Stein}).
Here the point is that one can attempt to raise levels
at primes that are ``Steinberg'' for $\rho$ and $\rho'$ and
Steinberg lifts have the property that the image of inertia
is unipotent and hence infinite. Thus Lemma \ref{artin} on which
Lemma \ref{conductors} relies does not apply.
\end{rmk}

\begin{rmk} A mod $pq$ descent result of the type suggested in
Remark 1 is not true for 2-dimensional
representations even in the simplest case when $K/\Q$ is a
quadratic extension. The principle for constructing a
counterexample is the following: Consider representations $\rho$
and $\rho'$ as above such that $\rho$ is ramified at a prime
$\ell \neq p,q$ and $>2$, and $\ell$ is not congruent to $\pm 1$
either mod $p$ or mod $q$, with the image of inertia
$\rho(I_{\ell})$ non-trivial and unipotent, $\rho'$ unramified at
$\ell$ and such that $\rho'({\rm Frob}_{\ell})$ has eigenvalues
with ratio $-\ell$. Let $\rho_{\ell}:=\rho|_{D_{\ell}}$ and
$\rho'_{\ell}:=\rho|_{D_{\ell}}$. Then it is easy from local
computations at $\ell$
that use structure of tame inertia to see that
while $\rho_{\ell}$ and $\rho'_{\ell}$ cannot arise from a single
Weil-Deligne representation $(\tau,N)$ as above, the restriction
of $\rho_{\ell}$ and $\rho'_{\ell}$ to the absolute Galois group
of the unramified quadratic extension of $\Q_{\ell}$ can arise in this
way. From this it is easy to construct counterexamples
remarking that if $\rho$ and $\rho'$ are irreducible on
restriction to  a quadratic extension $K$ of $\Q$ (with $\ell$ inert in $K$),
the extensions of $\rho|_{G_K}$ and $\rho'|_{G_K}$ to $G_\Q$ are
unique up to twisting by characters, and the condition that
$\rho'({\rm Frob}_{\ell})$ has eigenvalues with ratio $-\ell$ is
invariant under twisting. (The case of solvable mod $p$
automorphic descent for 2-dimensional representations is the main
theorem of \cite{K-MRL}.)
\end{rmk}

\begin{rmk} In a recent preprint of Manoharmayum \cite{JM} results towards Serre's conjecture mod 6 are proven.
\end{rmk}

\noindent {\bf Acknowledgement:} The authors wish to thank the referee for a number of remarks that helped improve the
presentation.

\end{document}